\documentclass{article}

\usepackage{arxiv}

\usepackage[utf8]{inputenc} 
\usepackage[T1]{fontenc}    
\usepackage{hyperref}       
\usepackage{url}            
\usepackage{booktabs}       
\usepackage{amsfonts}       
\usepackage{nicefrac}       
\usepackage{microtype}      
\usepackage{lipsum}
\usepackage{graphicx}
\usepackage{amsmath}
\usepackage{multirow}
\usepackage{titlesec}
\newtheorem{thrm}{Theorem}
\title{Analytical and Numerical Bifurcation Analysis of a Forest-Grassland Ecosystem Model with Human Interaction}

\author{
Konstantinos Spiliotis\\
Dipartimento di Agraria\\
Universit\'a degli Studi di Napoli Federico II\\
(Currently: Mathematical Institute\\
Rostock University,Germany)\\
\And
Lucia Russo\\
Consiglio Nazionale delle Ricerche\\
Napoli, Italy
\And
Francesco Giannino\\
Dipartimento di Agraria\\
Universita degli Studi di Napoli Federico II\\
Napoli, Italy
 \And
Constantinos Siettos\\
Dipartimento di Matematica e Applicazioni ``Renato Caccioppoli"\\
Universit\'a degli Studi di Napoli Federico II\\
Napoli, Italy\\
(Corresponding author: constantinos.siettos@unina.it)
}

\begin{document}
\maketitle

\begin{abstract} We perform both analytical and numerical bifurcation analysis of a forest-grassland ecosystem model coupled with human interaction. The model consists of two nonlinear ordinary differential equations incorporating the human perception of forest/grassland value. The system displays multiple steady states corresponding to different forest densities as well as regimes characterized by both stable and unstable limit cycles. We derive analytically the conditions with respect to the model parameters that give rise to various types of codimension-one criticalities such as transcritical, saddle-node, and Andronov-Hopf bifurcations and codimension-two criticalities such as cusp and Bogdanov-Takens bifurcations. We also perform a numerical continuation of the branches of limit cycles. By doing so, we reveal turning points of limit cycles marking the appearance/disappearance of sustained oscillations. These far-from-equilibrium criticalities that cannot be detected analytically give rise to the abrupt loss of the sustained oscillations, thus leading to another mechanism of catastrophic shifts.\end{abstract}

%

\section*{Introduction}	

Many ecosystems are characterized by a rich nonlinear behaviour including multistability as well the appearance of sustained oscillations and catastrophic shifts \cite{Scheffer_2001,Rietkerk_2004,Levin_1998,Jorgensen_2007,Capit_n_2010}. For example, catastrophic/regime shifts have been observed in shallow lakes\cite{Scheffer_2007}, the North Pacific, in the North Sea and Caribbean coral reefs \cite{deYoung2008}, in regions of northern Africa going abruptly from vegetated to desert conditions.\\ Frequently, the disturbances in ecosystems that drive such regime shifts are caused by human interactions \cite{Ling_2009,Genkai_Kato_2006,BURNEY_2005,Gordon_2008}. Thus, the understanding from a mathematical point of view of how disturbances/and or state changes influence the dynamics of ecosystems is very important in order to prevent and control ecological disasters \cite{Scheffer_2001,Villa_Mart_n_2015}.
From a nonlinear dynamics perspective, a sudden shift from one regime to another can be observed if the system has different coexisting stable regimes \cite{Scheffer_2001,Innes_2013,Villa_Mart_n_2015,Henderson_2016,Russo_2019}. In such cases, as a consequence of a vector state perturbation, the system may jump to a different stable state due to the crossing of its basin of attraction boundary \cite{Villa_Mart_n_2015}. Another way, in which sudden transitions may be observed is by parameter perturbations beyond a bifurcation point \cite{Villa_Mart_n_2015}. Thus, many studies in the field pinpoint the significance of bifurcation theory to explain, analyse and ultimately forecast beforehand sudden regime shifts \cite{Kooi2003,Scheffer_2001,Scheffer_2010,Bauch_2016,K_fi_2012,Ling_2009,Henderson_2016,Engler_2017,Cimatoribus_2013,Dijkstra_2009,Troost_2007}.
However, despite the well-established tools of bifurcation theory, many studies still use simple temporal simulations and/or simple linear numerical analysis to study models that are used to approximate real-world ecosystems. Recently, we have stressed the importance of using the arsenal of analytical and numerical methods of  bifurcation theory and particularly that of codimension-2 bifurcation analysis to systematically identify and characterize criticalities in ecological models and in general in physicalm biological and natural systems \cite{Russo_2019}. 
However, only few studies have used both analytical and numerical bifurcation analysis (i.e. not just temporal numerical simulations for the verification of the analytical findings, which do not suffice to identify and trace in the parameter space branches of far-from-equilibrium bifurcation points such as limit points of limit cycles), thus to study in a combined and complete way not only bifurcations at equilibrium but also bifurcations far-from-equilibrium (see for example \cite{Fujii_1982,Kooi2003,Troost_2007,Hamzah_2007,Dai_2010,Yu_2014}).

Here, we provide both analytical and numerical bifurcation results for a forest-grassland model with human interaction proposed by Innes et al. \cite{Innes_2013}. Forest-grassland mosaic ecosystems are typical examples of ecosystems where two species (forest and grass) compete for the same food (soil, sunlight and space) and their dynamics are strongly affected by human interactions \cite{Innes_2013,Bauch_2016}.
Despite the simplicity of the model, as we show, its nonlinear dynamical behavior is very rich including various types of codimension-one and codimension-two criticalities. We derive analytically the conditions for the onset of these criticalities including the codimension-one transcritical, saddle-node and Andronov-Hopf bifurcations as well as the codimension-two cusp and Bogdanov-Takens bifurcations. We also perform a two-parameter numerical bifurcation analysis to identify far-from-equilibrium criticalities such as limit points of limit cycles that cannot in principle be found analytically. By doing so, we provide a full picture of the mechanisms pertaining to the onset of codimension-one and codimension-two criticalities in a combined way. In the next section, we briefly present the mathematical model. We then provide the analytical results and consequently perform a numerical continuation of the limit cycles and construct the  two-parameter bifurcation diagram to fully characterize the different dynamical regimes. 


\section{Model Description}

The model describes the dynamics of a simple mosaic ecosystem of two species, grassland and forest under the human impact. Let $f$ and $g$ be the proportion of forest and grassland in the ecosystem, with $f+g=1$. Humans influence the fate of the evolution by selecting one of the species based on their abundance or lack. Let $x$ be the percentage of the human population that prefer forest; thus, $1-x$ is the ratio of people that prefer grassland. Each individual can alter his/her preference to forest or grassland according to some social-depended rules (educational or mimetic), or according to  the perceived value of the forest. The coupled model of ordinary differential equations (ODEs) reads \cite{Innes_2013}:
\begin{eqnarray}
\left\{ \begin{array}{l}
\frac{{df}}{{dt}} = w\left( f \right)\left( {1 - f} \right)f - \nu f - J\left( x \right)\\
\frac{{dx}}{{dt}} = sx\left( {1 - x} \right)U\left(f \right)
\end{array} \right.
\label{system}
\end{eqnarray}
where $\nu$ express the transition rate from forest to grassland due to natural processes and $w(f)$ is the rate of transition from grassland to forest mainly induced due to fire modelled by the following relation:
\begin{equation}
w\left( f \right) = \frac{c}{{1 + {e^{ - k\frac{f}{{1 - f}} + b}}}}
\label{fire}
\end{equation}
The $w$ function is of sigmoid type while the properties of its shape depend from the parameters $b$, $c$ and $k$; $c$ gives the maximum value, while $b$, $k$ are related to the rate of activation of $w$.

The human influence is represented in the above system of ODEs by the term $J(x)= h(1-2x)$, where $h$ is the potential magnitude of human influence on the ecosystem. The human interaction has a positive or negative feedback on the perceived value of the forest: if the majority prefers forest (i.e. $x>1/2$) the society will enhance forest recreation (i.e $-J(x)>0$) by reforestation; instead if social opinion favours grassland then the human impact leads to deforestation (i.e $-J(x)<0$). The perceived value of the forest $U(f)$ is given by:  $U(f)=(1-f)-f=1-2f$.\\
For more details on the model see Innes et al. 2013\cite{Innes_2013}.

\section{Analytical Results} \label{Model}
We first give a geometrical description for the multiplicity of fixed point solutions. Then, we derive analytically the fixed point solutions and characterize their stability. Consequently, we proceed with the construction of the one-parameter bifurcation diagrams with respect to the human influence $h$ as well as the $k$ parameter. Then, we provide the conditions for the existence of  the saddle-node and Andronov-Hopf bifurcations. 
 
 
 \subsection{Fixed Point Solutions-Stability Analysis}
Setting the derivative of the second equation of Eq.(\ref{system}) equal to zero, we get three solutions:  $x=0$, the $x=1$ and $f=1/2$.\\

Starting with $x=0$ and substituting to the first equation of Eq.(\ref{system}) and setting this to zero, we obtain 
\begin{equation}
w\left( f \right)\left( {1 - f} \right)f - \nu f - h=0
\end{equation}
The above equation defines a nonlinear function of fixed points for the forest density $f$ ($0\leq f \leq 1$) with respect to $h$  when $x=0$:
\begin{equation}
h=h\left( f\right)=w\left( f \right)\left( {1 - f} \right)f - \nu f
\label{x_0}
\end{equation}

\begin{figure}[pb]
\centerline{\includegraphics[width=3 in]{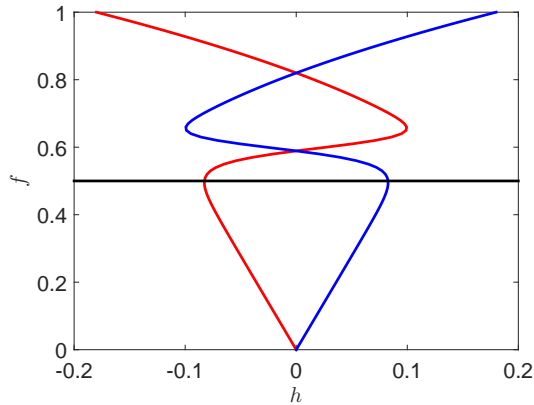}}
\vspace*{8pt}
\caption{Three branches of fixed points solutions, $(f(h), x=0)$ (marked with red), $(f(h), x=1)$ (marked with blue) and the horizontal line $(f=1/2$,$x=0$, (or $x=1$), (marked with black). We used $c=1$, $s=10$, $b=11$, $k=6.5$ and $\nu=0.2$. The branches are constructed according to (\ref{x_0}) and (\ref{x_1}).}
\label{figx_0_x_1}
\end{figure}

In Fig.\ref{figx_0_x_1} we have plot $h(f)$ for a specific choice of the values of the other parameters of the model.

For $x=1$, we get a symmetric expression for the parameter $h$ for the second branch of fixed points $(f(h),x=1)$ i.e.
\begin{equation}
h=h\left( f\right)=-w\left( f \right)\left( {1 - f} \right)f + \nu f
\label{x_1}
\end{equation}

This is the mirror case of the fixed points corresponding to $x=0$ with respect to the y-axis (see blue line in Fig.\ref{figx_0_x_1}).\\

Finally, for $f=1/2$ and by setting zero the first equation of the system ~(\ref{system}), we get:
\begin{equation}
w\left( 1/2\right)\left( {1/4} \right) - \nu/2 =h\left( {1-2x } \right) \Leftrightarrow  \frac{c}{{1 + {e^{ b-k}}}}-2\nu=4h(1-2x)
\label{hrel}
\end{equation}

The above gives rise to another family of solutions of fixed points $(f=1/2,x)$ with $x$ given by:
\begin{equation}
x=1/2\left(1- \frac{1}{4h}\left( \frac{c}{{1 + {e^{ b-k}}}}-2\nu\right)\right)
\label{f_12x}
\end{equation}

\subsubsection{A Geometrical Description of the Multiplicity of Fixed Points}

\begin{figure}[ht]
\centerline{\includegraphics[width=6 in]{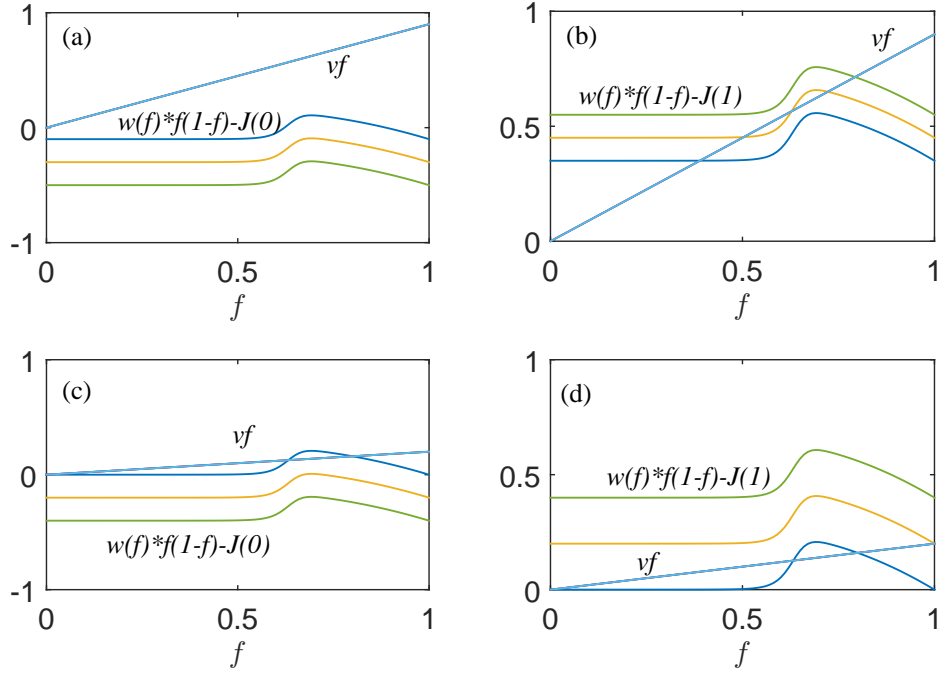}}
\vspace*{8pt}
\caption{Fixed point solutions for $s=10$, $b=11$, $c=1$, $k=6.5$, (a) $x=0$ and $\nu=0.9$. For this choice of parameter values, there are no fixed points solutions for positive values of $h$.
(b) $x=1$. Depending on $h$, the number of fixed points varies from one to three. (c) Fixed point solutions for $x=0$ and $\nu=0.2$. The parameter $h$ varies from zero to three. Depending on the value of $h$, the number of fixed points varies from zero to three. (d) $x=1,\nu=0.2$.}
\label{fixedpoints}
\end{figure}

Fig.\ref{fixedpoints} shows the fixed point solutions for $x=0$, $x=1$ as they result from  Eq.(\ref{x_0}) and Eq.(\ref{x_1}) as intersection points of the graphs of the functions $\nu f$ and $w\left( f \right)\left( {1 - f} \right)f-J(x)$ for various choices of the values of the other parameters. Fig.\ref{fixedpoints}(a) depicts the intersection of the graphs for $k=6.5$, $b=11$, $c=1$, $\nu=0.9$). Clearly, for $h>0$ there are not fixed points for the $x=0$ branch of solutions. Fig.\ref{fixedpoints}(b) depicts the intersection of the graphs for the $x=1$ branch of solutions. As $h$ increases, the number of steady state regimes first increases from one to three and then decreases to zero. Finally, it should be considered that the branch of steady state solutions $f=1/2$ exists for all values of $h$. 

For $k=6.5$, $\nu=0.2$ the intersections of the graphs are depicted in Fig.\ref{fixedpoints}(c),(d) as $h$ is varied. In particular, in Fig.\ref{fixedpoints}(c) we show fixed points  corresponding to $x=0$ and in Fig.\ref{fixedpoints}(d) the ones corresponding to $x=1$. As it is shown, for $x=0$, as $h$ is larger than $0$, the number of solutions goes from three to two, indicating that $h=0$ is a bifurcation point. When $h>0.07$, there are no solutions corresponding to $x=0$. For the case of $x=1$, when $h>0$, there are three solutions, two of them disappearing at $h=0.1$. Again for very high values of $h$, there are no solutions corresponding to $x=1$.\\ 
While this graphical analysis can give a quick idea of the number of solutions for different values of $h$, it cannot give a precise information about the stability and type of bifurcations that may appear. Moreover, dynamic regimes like periodic solutions cannot be detected. Thus, in what follows, we perform a systematic analytical study of the stability of fixed point solutions and their bifurcations.
 
\subsubsection{Stability Analysis of Fixed Points}
The stability of the above families of fixed points is determined from the eigenvalues $\lambda_1, \lambda_2$ of the Jacobian matrix of~\ref{system} evaluated at the fixed points.\\
The Jacobian matrix of~\ref{system} reads:
\begin{equation}
J=
    \quad
\begin{pmatrix} 
w'\left(f\right) \cdot(f-f^{2})+w\left(f\right)\cdot(1-2f)-\nu & 2h \\
-2sx(1-x) & s(1-2x)\cdot \left(1-2f\right) 
\end{pmatrix}
\label{Jac}
\end{equation}

Clearly, when $x=0$ or $x=1$, the Jacobian matrix has an upper triangular form since $J_{2,1}$ is vanished. In these cases, the eigenvalues  of the Jacobian matrix are real and equal to the diagonal elements. This implies that any loss of hyperbolicity (for $x=0$ or $x=1$) cannot be due to Andronov-Hopf bifurcations.\\
The determinant of the Jacobian matrix for $x=0$ is given by:
\begin{equation}
    det(J)|_{x=0}= 
    \begin{vmatrix} 
    w'\left(f\right) \cdot(f-f^{2})+w\left(f\right)\cdot(1-2f)-\nu & 2h \\
0 & s\cdot \left(1-2f\right) 
    \end{vmatrix}\\
    \label{Jac_x0}
\end{equation}

or

\begin{equation}
     det(J)|_{x=0} =\left( w'\left(f\right)\cdot(f-f^{2})+w\left(f\right)\cdot(1-2f)-v\right)\cdot s\left( 1-2f \right)
     \label{Jac_x0b}
\end{equation}

Note that:

\begin{equation}
    det(J)|_{x=1}=-det(J)|_{x=0}
    \label{opos_stab}
\end{equation}

Eq.(\ref{opos_stab}) provides the stability relation between the $x=0$ and $x=1$ branches: the stability of one branch implies the instability of the other. Thus, it suffices to study the stability properties of just one of them, for example the $x=0$ solution branch. From Eq.(\ref{fire}), the derivative of $w$ with respect to $f$ is given by: 
\begin{equation}
w'\left(f\right) = \frac{c\cdot k e^{ - k\frac{f}{{1 - f}} + b}}{\left({1 + {e^{ - k\frac{f}{{1 - f}} + b}}}\right)^2}\cdot \frac{1}{\left(1-f\right)^2}>0
\end{equation}

Note that both $w(f)$ and $w'\left(f\right)$ are always positive (i.e. $w$ is a positive increasing function of $f$). 
The trace of $J$ at $x=0$ is given by:

\begin{equation}
     Tr(J)|_{x=0}= \left( w'\left(f\right)\cdot(f-f^{2})+w\left(f\right)\cdot(1-2f)-v\right)+ s\left( 1-2f \right)
\end{equation}

Thus, for $x=0$, we examine separately two cases that mark the onset of a critical point: (a) $J_{2,2}=0$  (resulting from $f=1/2$ (see Eq.(\ref{Jac_x0}))) and (b) $J_{1,1}=0$. \\
As we prove in the next section, the case (a) corresponds to the appearance of a transcritical bifurcation, while the case (b) corresponds to the appearance of a saddle-node bifurcation. 

\subsection{Transcritical Bifurcations}
We examine the case when the two branches of fixed points intersect at $(f^*,x^*)=(1/2,0)$. Note by Eq.(\ref{f_12x}) and Eq.(\ref{hrel}), that the value of the parameter $h$ that results to $x=0$ in Eq.(\ref{f_12x}) is given by: $h^*=w(1/2)/4-\nu/2$ or equivalently by:


\begin{equation}
    h^*=\frac{c}{4(1+e^{b-k})}-\nu /2
    \label{h*}
\end{equation}

That is, for the above parameter value $h^*$, the above families of branches of fixed points (i.e. the ones defined by (a) $x=0$ and $h(f)$ (given by Eq.(\ref{x_0})) and  (b) $f=1/2$ and $x$ (given by Eq.(\ref{f_12x})) intersect at $(f=1/2,x=0)$.

Depending on the values of the other parameters, i.e. $c$, $b$, $k$, the value of $h^*$ may change sign.\\
The upper triangular form of the Jacobian matrix evaluated at $(f^*,x^*)=(1/2,0)$ implies that $J_{22}=\lambda_2=0$, while the second eigenvalue is given by (see Eq.(\ref{Jac_x0})):
$\lambda_1= w'\left(1/2\right)\cdot(1/4)-\nu=\frac{c\cdot k e^{b - k}}{\left({1 + {e^{ b- k}}}\right)^2}-\nu$.\\
Clearly if $\lambda_1$ is of constant sign around $f=1/2$:

\begin{equation}
    \frac{c\cdot k e^{b - k}}{\left({1 + {e^{ b- k}}}\right)^2}-\nu <>0
    \label{casNotEq}
\end{equation}

then the second eigenvalue $J_{22}=\lambda_2=s(1-2f)$ changes sign with respect to $f$, when $f\in (\frac{1}{2}-\epsilon,\frac{1}{2}+\epsilon)$: (when $\lambda_2=s(1-2f)>0\Leftrightarrow f<1/2$ and when $\lambda_2=s(1-2f)<0\Leftrightarrow f>1/2$).\\

If $\lambda_1<0$,i.e. when
\begin{equation}
\lambda_1=w'\left(1/2\right)\cdot(1/4)-\nu<0 \Leftrightarrow \frac{c\cdot k e^{b - k}}{\left({1 + {e^{ b- k}}}\right)^2}-\nu < 0
    \label{casNeg}
\end{equation}
then we have a change of stability of solutions 
from unstable $(\lambda_2>0 (f<1/2), \lambda_1<0)$ to stable $(\lambda_2<0  
 (f>1/2), \lambda_1<0)$.\\

If $\lambda_1>0$, i.e. when
\begin{equation}
   \lambda_1= w'\left(1/2\right)\cdot(1/4)-\nu>0 \Leftrightarrow \frac{c\cdot k e^{b - k}}{\left({1 + {e^{ b- k}}}\right)^2}-\nu > 0
     \label{casPos}
\end{equation}
the transition is from unstable nodes ($f<1/2$) (two positives eigenvalues) to saddle nodes (two eigenvalues with opposite sign).\\
We summarise the above results in the following Theorem. 
\begin{thrm}
The point $(f^*,x^*)=(1/2,0)$ is a critical point for the branch of fixed point solutions defined by $x=0$. The value of the parameter $h^*$ that defines this criticality is given by $h^*=\frac{c}{4(1+e^{b-k})}-\nu /2$. In the case where (\ref{casNeg}) holds true, there is a change of stability from unstable states when $f<1/2$ to stable states when $f>1/2$. If Eq.(\ref{casPos}) holds, then the $x=0$ branch remains locally unstable.
\end{thrm}

To determine the type of bifurcation at $(f^*,x^*)=(1/2,0)$, we proceed as follows.

The Jacobian matrix at $(f^*,x^*)=(1/2,0)$ reads:

\begin{equation}
J=
    \quad
\begin{pmatrix} 
w'\left(1/2\right) 1/4-\nu &  \hspace{0.5cm} 2h \\
0 & \hspace{0.5cm}0 
\end{pmatrix}
\label{Jax121}
\end{equation}
meaning that $rank(J)_{f=1/2,x=0}=1$ and the derivative of the system Eq.(\ref{system}) with respect to the parameter $h$ is $ J_h(1/2,0)=e_{1}=[1,0]^T$, meaning that 
\begin{equation}
    J_h(1/2,0)\in Range(J(1/2,0))=span \{[1,0]^T \}.
    \label{condTran}
\end{equation}

Then, $(f^*,x^*)=(1/2,0)$ is a simple bifurcation point \cite{Seydel_2010} and since there are exactly two intersecting branches, we conclude that this type of bifurcation is a transcritical bifurcation.\\
The stability of the $f=1/2$ branch can again be determined with the aid of Eq.(\ref{f_12x}) and Eq.(\ref{h*}). Thus, the solution $x$ for the branch $f=1/2$ can be re-written as:
\begin{equation}
x=1/2 (1- \frac{h^*}{h})
\label{f_12xa}
\end{equation}

For $h<h^*<0\implies1>h^*/h \implies x>0$, (while if $0>h>h^*$, we have the opposite (non-physical) scenario i.e. $x<0$). The determinant along the $f=1/2$ branch of solutions is:
\begin{equation}
    det(J)|_{f=1/2}= 4s\cdot h\cdot x(1-x)
\end{equation}

Note that $1-x>0$ and assuming that $s$ is constant, if $h<h^*<0 \implies det(J)|_{f=1/2}<0$, and, if $h^*<h<0 \implies det(J)|_{f=1/2}>0$. \\

Here, we comment that the stability properties of this branch are not independent from the stability of the $x=0$ branch. The trace of the Jacobian at $f=1/2$ is:
\begin{equation}
     Tr(J)|_{f=1/2}=  w'\left(1/2\right)\cdot1/4-\nu
\end{equation}
Then, assuming that $Tr(J)|_{f=1/2}<0$ (and in combination with the sign of the determinant) this is translated as a change from unstable $(h<h^*<0)$ to stable $(h^*<h<0)$ fixed point solutions. But the assumption $Tr(J)|_{f=1/2}<0$ is identical with Eq.(\ref{casNeg}), which defines the change of stability on the $x=0$ branch around the critical point $(f^*,x^*)=(1/2,0)$, from unstable $(f<1/2)$ to stable $(f>1/2)$ solutions (see Eq.(\ref{casNeg}-\ref{casPos})).\\
To conclude, at the fixed point $(x^*,f^*)=(0,1/2)$ we have an exchange of stability, which is characteristic of transcritical bifurcation. We state the following theorem:

\begin{thrm}
The system undergoes a transcritical bifurcation at $(f^*,x^*)=(1/2,0)$. In the cases where Eq.(\ref{casNeg}) and $h^*<0$ hold true, then, as $h$ and $f$ increase, the two branches $x=0$ and $f=1/2$ exchange their stability from unstable to stable state and vice versa.  
\end{thrm}


\subsection{Saddle-Node and Cusp Bifurcations}

\subsubsection{Bifurcations on the $x=0$ branch of fixed point solutions}
Here, we examine the second case ($J_{11}=0$) for the existence of non-hyperbolic points on the branch $x=0$, (see Eq.(\ref{Jac_x0})).\\
For this case, we obtain:
\begin{equation}
   J_{11}(f)=w'(f)\cdot (f-f^2)+w(f)\cdot (1-2f)-\nu=0
    \label{J110a}
\end{equation}
or equivalently
\begin{equation}
   J_{11}(f)=\frac{c\cdot k e^{ - k\frac{f}{{1 - f}} + b}}{\left({1 + {e^{ - k\frac{f}{{1 - f}} + b}}}\right)^2}\cdot \frac{f}{\left(1-f\right)}+\frac{c}{{1 + {e^{ - k\frac{f}{{1 - f}} + b}}}}\cdot (1-2f)-\nu=0
    \label{J110b}
\end{equation}
\begin{figure}[pb]
\centerline{\includegraphics[width=3.5 in]{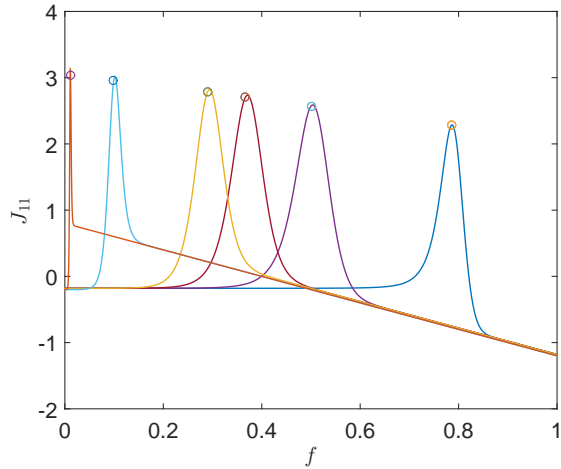}}
\vspace*{8pt}
\caption{Graphical representation of $J_{11}$ function for $k$=3, 11, 19, 27, 100, 1000 (from right to left). with circle marker we depict the approximate value of $max(j_{11})$ resulting from Eq.(\ref{Appr_max}).}
\label{J11}
\end{figure}

The system of Eq.(\ref{x_0}),Eq.(\ref{J110b}) defines the appearance of critical points (the Jacobian becomes singular). The two nonlinear equations can be solved numerically (e.g. by the Newton-Raphson method) with respect to the unknowns $h$, $f$. We will show analytically, that depending on the value of the parameters $h$, $k$ and $\nu>0$, the above nonlinear system can have at most two solutions. Furthermore, we show that the corresponding criticalities are saddle-node bifurcations.\\
From Eq.(\ref{J110b}) we see that as $f\rightarrow 0^+$, $J_{11}\rightarrow \frac{c}{1+e^{b}} -\nu \approx -\nu<0$. As $f\rightarrow 1^-$ then $y=\frac{f}{1-f}\rightarrow +\infty$, then $ e^{ - k\frac{f}{{1 - f}} + b}\cdot \frac{f}{1-f}=e^{ - k y + b}\cdot y=\frac{y}{e^{ k y - b}}\rightarrow 0$, concluding $J_{11}\rightarrow (1-2f)-\nu=-1-\nu<0$.
Fig.\ref{J11} depicts $J_{11}$ for different values of the parameter $k$, which confirms the asymptotical behavior of $J_{11}$ at the boundaries $f=0$ and $f=1$. In the neighborhood of $f=0$, it converges to a constant function the $-\nu$, while, in the neighborhood of $f=1$, it converges asymptotically to $(1-2f)-\nu$. Between the two boundaries, $J_{1,1}$ has a bell shape resulting from the first part of $J_{11}$, which is $w'(f) \frac{f}{\left(1-f\right)}$. Thus, $J_{11}$ has two monotonic intervals: the $[0, f^*]$ in which it increases and the $[f^*, 1]$ in which decreases. Both $w'$ and $\frac{f}{1-f}$ are positive functions and $\frac{f}{1-f}$ increases. We study the monotonicity properties of the first part of $w'$. We consider the function $L(f)$:           
\begin{equation}
   L(f)=\frac{ e^{ - k\frac{f}{{1 - f}} + b}}{\left({1 + {e^{ - k\frac{f}{{1 - f}} + b}}}\right)^2}
   \label{Lf}
\end{equation}
then by setting $z=e^{-k\frac{f}{\left(1-f\right)}+b}$, Eq.(\ref{Lf}) reads
\begin{equation}
   L(z)=\frac{ z}{(1 + z)^2},
    \label{Lz}
\end{equation}
where now $z\in[0, e^b]$. Then, we get the location of the maximum at $z^*=1$ or equivalently $\frac{\tilde{f^*}}{1-\tilde{f^*}} = \frac{b}{k} \Leftrightarrow \tilde{f^*}=\frac{b}{b+k}$. Substituting this in Eq.(\ref{J110b}) we obtain:
\begin{equation}
    J_{1,1}(\tilde{f^*})=\tilde{J}_{1,1}=\frac{c b}{4}+\frac{k-b}{2(k+b)}-\nu.
    \label{Appr_max}
\end{equation}
Fig.\ref{J11} shows the $\tilde{J}_{11}$ (marked with circles) for different values of the parameter $k$. Note that the values of $J_{11}(\tilde{f^*})$ are in a very good agreement with actual maxima of $J_{11}$. Depending on the values of $\nu$, $k$ the maximum of $J_{11}$ can be positive:
 \begin{equation}
    max(J_{11})>0.
    \label{maxJpos}
\end{equation}
Then from the  mean value theorem for continuous functions, this implies that the function $J_{11}$ has two roots correspond to turning points. Since
 \begin{equation}
     max(J_{11})\geq J_{11}(\tilde{f^*})=\frac{c b}{4}+\frac{k-b}{2(k+b)}-\nu \geq \frac{cb}{4}-\frac{1}{2}-\nu,
     \label{maxJposcond}
 \end{equation}
then in order to have two turning points, a sufficient condition for the parameters $b$, $k$, $\nu$ is given by:
 \begin{equation}
     \frac{cb}{4}-\frac{1}{2}-\nu>0.
 \end{equation}
 
In the critical condition where:
\begin{equation*}
    max(J_{11})=0
\end{equation*}

or approximately:

 \begin{equation}
     \frac{c b}{4}+\frac{k-b}{2(k+b)}=\nu,
     \label{cuspcond}
 \end{equation}
 
the two turning points collapse to a cusp bifurcation.\\
In the first case, where the function $J_{11}$ has two roots, say $\rho_{1,2}$, the augmented matrix $[J,J_h]$ calculated on the roots of $J_{11}$ reads:

 \begin{equation}
    J_{augm}=
    \quad
\begin{pmatrix} 
0 &  2h & -1\\
0 & \hspace{0.5cm}s (1-2\rho_{1,2}) & \hspace{0.3cm} 0 \\
\end{pmatrix} 
 \end{equation}
 
and for $\rho_{1,2}\neq1/2$, $J_{augm}$ has a rank $r(J_{augm})=2$ implying that the the roots of $J_{11}$ are turning points \cite{Seydel_2010}. The stability properties of the $x=0$ branch near the turning points results in a straightforward manner from the signs of eigenvalues. We assume that the roots of $J_{11}$ have the following order $1/2<\rho_1<\rho_2$. This ordering can be achieved if Eq.(\ref{casNeg}) holds true and simultaneously $\frac{b}{k+b}>1/2 \Leftrightarrow b>k $ (which for example is the case with $b=11, k=6.5$). Table \ref{tablesignx0} summarizes the stability along the $x=0$ with respect to the signs of the eigenvalues.  
\begin{table}[ht]
\caption{Stability analysis of the $x=0$ branch in the case where the Jacobian becomes singular at $f=1/2,\rho_1,\rho_2$, where the roots have the following order: $1/2<\rho_1<\rho_2$}
    \centering
 \begin{tabular}{ |c|c|c|c|c|}
 \hline
 \multicolumn{5}{|c|}{f \hspace{2.1cm}1/2 \hspace{0.7cm} $\rho_{1}$  \hspace{1cm} $\rho_{2}$ \hspace{1.cm}} \\
 \hline
 $\lambda_1=J_{11}$&-&-&+&-\\
 \hline
 $\lambda_2=s(1-2f)$ &+&-&-&-\\
 \hline
 \hline
 Stability & unstable & stable & unstable & stable \\
 \hline
\end{tabular}
\label{tablesignx0}
\end{table}

We summarise the above analytical results in the following theorems. 
\begin{thrm}
Assuming that Eq.(\ref{casNotEq}) and  Eq.(\ref{maxJpos}) hold true, then the $x=0$ branch of fixed points posses two turning points. In addition, if the conditions $b>k$ and Eq.(\ref{casNeg}) are imposed i.e. ($1/2<\rho_1<\rho_2$), then stability properties are determined according to Table \ref{tablesignx0}.
\end{thrm}

\begin{thrm}
Assuming that the condition (\ref{cuspcond}) holds for the model parameters, then the $x=0$ branch of fixed points posses a cusp bifurcation.
\end{thrm}

\subsubsection{Bifurcations on the $x=1$ branch of fixed point solutions}

In this section, we conduct our analysis for the second symmetric branch of fixed point solutions, $x=1$. Following the previous analysis, we study the transcritical bifurcation at the point $(1/2,1)$ and the stability properties of turning points. 

The determinant of Jacobian matrix along the $x=1$ is given by:
\begin{equation}
    det(J)|_{x=1}= 
    \begin{vmatrix} 
    w'\left(f\right) \cdot(f-f^{2})+w\left(f\right)\cdot(1-2f)-\nu & 2h \\
0 & -s\cdot \left(1-2f\right) 
    \end{vmatrix}\\
    \label{Jac_x1}
\end{equation}

and similar to the $x=0$ case, the triangular form of the Jacobian matrix gives immediately the eigenvalues on the branch: $\lambda_1' =J_{11}=\lambda_1$ and  $\lambda'_2=J_{22}=-s(1-2f)=-\lambda_2$. For the simplicity of the presentation,  we suppose that hold true the same assumptions for $\lambda_1$ (i.e. Eq.(\ref{casNeg})) and for the ordering of the roots of $J_{11}$  as in the previous case.\\
The determinant becomes singular at $(f, x)=(1/2,1)$, while if $f>1/2 \implies \lambda'_2>0$  and if $f<1/2 \implies \lambda'_2<0$. Near the point $(1/2,1)$ and as the parameter $h$ increases, the $x=1$ branch changes stability from stable to unstable. The critical value of the parameter $h$ on the bifurcation point resulting from Eq.(\ref{x_1}) by setting $f=1/2$ is given by:  $h^{**}=-h^*=-w(1/2)/4+\nu/2$ or

\begin{equation}
    h^{**}=-\frac{c}{4(1+e^{b-k})}+\nu /2
    \label{h**}
\end{equation}

The Jacobian matrix evaluated at the point $(1/2,1)$ is exactly the same as the one evaluated at $x=0 $ (Eq.(\ref{Jax121})) and the condition Eq.(\ref{condTran}) continues to be valid, concluding that at $(1/2,1)$ a transcritical bifurcation appears. 

The stability of the $f=1/2$ branch near the critical point $(1/2,1)$ can be determined as follows. By substituting $h^{**}$ (given by Eq.(\ref{h**})) into Eq.(\ref{f_12x}) we get:

\begin{equation}
x=1/2 (1+ \frac{h^{**}}{h})
\label{f_12xb}
\end{equation}

Assuming that $h^{**}>0$, then for $h>h^{**}>0\implies1>h^{**}/h \implies 1-x=1/2(1-h^{**}/h) >0$. The determinant reads:
\begin{equation}
    det(J)|_{f=1/2}= 4s\cdot h\cdot x(1-x)
\end{equation}

concluding that $det(J)|_{f=1/2}>0 $ for $h>h^{**}>0$ and $det(J)|_{f=1/2}<0 $ for $h^{**}>h>0$. Since the condition given by Eq.(\ref{casNeg}) holds true, the sign of the trace of the Jacobian reads:

\begin{equation}
     Tr(J)|_{f=1/2}=  w'\left(1/2\right)\cdot1/4-\nu<0.
\end{equation}

The above implies that as $h$ increases, the stability of the branch $f=1/2$ changes from  unstable $(f<1/2)$ to stable $(f>1/2)$.


We summarize the above results with the following theorem.

\begin{thrm}
The system (~\ref{Model}) undergoes a transcritical bifurcation at $(f^*,x^*)=(1/2,1)$. In the cases where Eq.(\ref{casNeg}) and $h^{**}>0$ hold true, then, the two branches $x=1$ and $f=1/2$ exchange stability at that bifurcation point. The value of the parameter $h$ on the bifurcation point is given by $h^{**}=-\frac{c}{4(1+e^{b-k})}+\nu /2$.
\end{thrm}

The stability properties of the $x=1$ branch can easily be calculated from Eq.(\ref{opos_stab}) and from the analysis of the $J_{11}$ function. Table \ref{tablesignx1} summarizes the stability along the branch $x=1$ with respect to the sign of the eigenvalues. We observe that for $f>1/2$, the branch of solutions remains unstable. At the turning points $\rho_{1,2}$, the transition is from saddle to unstable nodes (in $\rho_1$) and again back to saddles.

\begin{thrm}
Assuming that Eq.(\ref{casNotEq}) and Eq.(\ref{maxJpos}) hold, then the $x=1$ branch of solutions posses two turning points. In addition if the conditions $b>k$ and Eq.(\ref{casNeg}) are imposed i.e.($1/2<\rho_1<\rho_2$), then the stability properties of the branch are determined according to Table \ref{tablesignx1}.
\end{thrm}

\begin{table}[ht]
\caption{Stability of the $x=1$ branch in the case where the Jacobian becomes singular at $f=1/2$, $\rho_1$, $\rho_2$ where the roots have the following order: $1/2<\rho_1<\rho_2$}
    \centering
 \begin{tabular}{ |c|c|c|c|c|}
 \hline
 \multicolumn{5}{|c|}{f \hspace{2.1cm}1/2 \hspace{0.7cm} $\rho_{1}$  \hspace{1cm} $\rho_{2}$ \hspace{1.cm}} \\
 \hline
 $\lambda'_1=J_{11}$&-&-&+&-\\
 \hline
 $\lambda'_2=s(1-2f)$ &-&+&+&+\\
 \hline
 \hline
 Stability & stable & unstable & unstable & unstable \\
 \hline
\end{tabular}
\label{tablesignx1}
\end{table}

\begin{figure}[ht]
\centerline{\includegraphics[width=3.5 in]{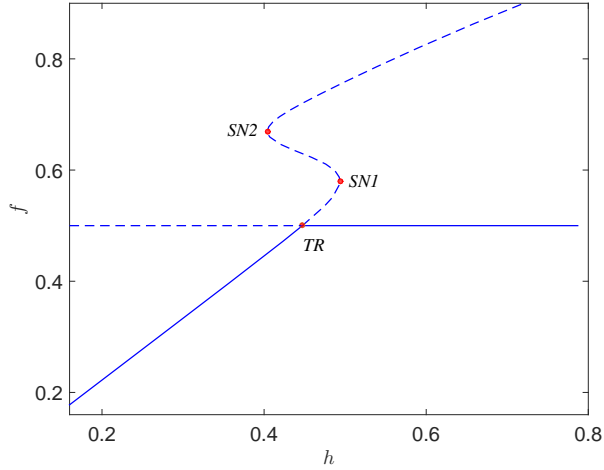}}
\vspace*{8pt}
\caption{Bifurcation diagram of $f$ with respect to the parameter $h$. The ``S"-type branch bifurcates via a transcritical bifurcation $TR1$ at the point $(f,x,h)=(1/2,1,0.447)$ and it exchanges the stability with the horizontal branch. The ``S" branch contains two more saddle node bifurcation points at $(f,x,h)=(0.58,1,0.494)$ ($SN1$) and at $(f,x,h)=(0.67,1,0.404)$($SN2$).}
\label{bifx1v09}
\end{figure}

As an example of the above analysis, we present the bifurcation diagram of $f$ with respect to $h$ for $k=6.5$, $\nu=0.9$, $s=10$ (Fig.\ref{bifx1v09}). Solid lines depict stable solutions, while dashed lines unstable ones. Note that for this choice of the values of $s$, $k$, $\nu$, there are no fixed point solutions for $x=0$ ($h>0$); hence this bifurcation branch corresponds to the locus of steady states with $x=1$. The horizontal line corresponds to the regime fixed points with $f=1/2$. Starting from small values of $h$, the system has two solutions, one stable corresponding to $x=1$ and another unstable for $f=1/2$. As $h$ increases, the ''S" branch  crosses the horizontal line and the two branches exchange their stability through a transcritical bifurcation ($TR$) at $(f,x,h)=(0.5,1,0.447)$. The ''S" branch also displays two saddle-node bifurcation points ($SN1$) at $(f,x,h)=(0.58,1,0.494)$ and ($SN2$) at $(f,x,h)=(0.67,1,0.404)$. Between ($SN1$) and ($SN2$), the system displays four steady states and just one of these is stable.

\begin{table}[ht]
\caption{Critical Values of parameter $h$ corresponding to the bifurcation diagram shown in Fig.\ref{bifx1v09}; the other values of the parameters were set as $k=6.5$, $\nu=0.9$.}
\begin{tabular}{|c|c|c|}
\hline
Bifurcation Point & Kind of Bifurcation & Branch \\
$(f^{*},x^{*},h^{*})$ &  & \\ 
\hline
$(0.5,1,0.447)$ & TR & $f=1/2, x=1$\\
\hline
(0.58,1,0.494) & SN1 & $x=1$\\
\hline
(0.67,1,0.404) & SN2 & $x=1$\\ 
\hline
\end{tabular}
\end{table}

\begin{figure}[ht]
\centerline{\includegraphics[width=5.5 in]{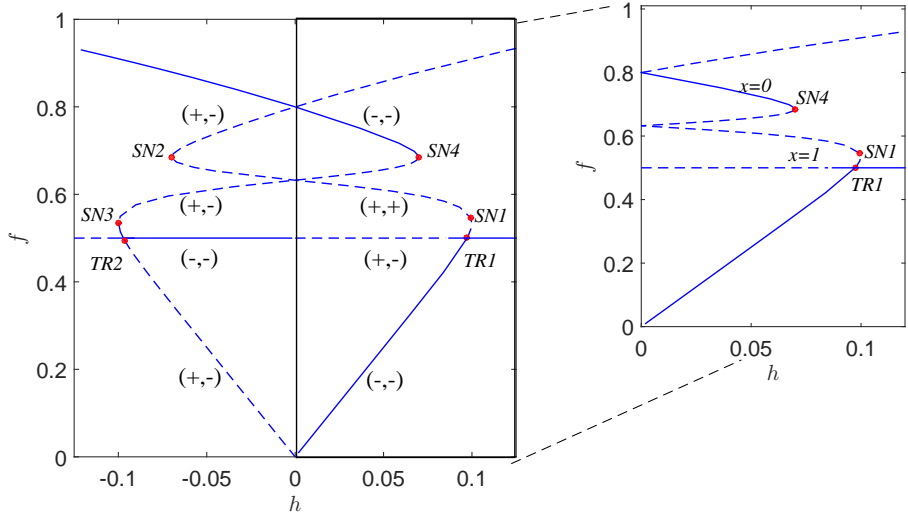}}
\vspace*{8pt}
\caption{Bifurcation diagram of $f$ with respect to the parameter $h$; the values of the other parameters are set as $k=6.5$, $b=11$, $c=1$, $s=10$. Solid lines correspond to stable fixed points and dashed lines to unstable fixed points. $SN$ denote saddle-node bifurcations and $TR$ transcritical bifurcations. The diagram consists of two branches one of ``S"-type ($x=1$ ) and one inverted ``S" ($x=0$). The exact values of bifurcations point are given in Table\ref{tab_bifot}. The inset on the right shows the bifurcations diagram only for positive values of $h$.}
\label{bifx0x1}
\end{figure}

Returning back to Fig.\ref{figx_0_x_1}, we are now in the position to characterize the types of bifurcations and stability of the solution branches. Fixing the values of the parameters at $s=10$, $k=6.5$, $b=11$, $\nu=0.2$, the annotated  bifurcation diagram of $f$ with respect to the parameter $h$ is shown in Fig.\ref{bifx0x1}.  

Solid lines depict stable fixed points, while dashed lines depict unstable ones. As discussed above, there are three branches of steady states: one inverted ``S branch corresponding to $x=0$, an ``S" branch  corresponding to the locus of steady states with $x=1$, and the horizontal line corresponding to the fixed points with $f=0.5$. 
In order to understand, from a nonlinear dynamical point of view, the behaviour of the $x=0$ branch of fixed points, we show the bifurcation diagram also for negative values of $h$. In particular, in Fig.\ref{bifx0x1} the stability of each fixed point is reported: $(-,-)$ are stable nodes with two negative eigenvalues; $(-,+)$ are saddles with one negative and one positive eigenvalue, and $(+,+)$ are unstable nodes with two positive eigenvalues. It is apparent, that the branch of the fixed points $x=0$ appears with a mirror ``S"-shape with respect to the $x=1$ ```S" branch (as we proved in the previous section). Both the ``S"-shaped branches ($x=0$ and $x=1$) cross the horizontal line $f=1/2$ giving rise to two transcritical bifurcations ($TR1$ for the $x=1$ branch, $TR2$ for the $x=0$ branch) on which the branches of fixed points exchange their stability.
Thus, looking just at the positive values of $h$, up to the value corresponding to the saddle node bifurcation $SN4$, there are six steady states: two stable nodes (one for $x=0$ and one for $x=1$), three saddles (one for $x=0$,one for $x=1$ and one for $f=1/2$) and one unstable node (for $x=0$). Multistability and multiplicity of fixed points is observed in a wide range of parameters, whereas as the parameter $h$ passes the value of one corresponding to the saddle-node point $SN1$, only one stable steady state exists. 

\begin{table}[ht]
\caption{Critical Values of $h$ parameter.}
{\begin{tabular}{|c|c|c|}
\hline
Bifurcation Point & Kind of Bifurcation & Branch \\
$(f^*,x^*,h^*)$ &  & \\ 
\hline
$(0.5,1,0.096)$ & TR1 & $f=1/2, x=1$\\
\hline
$(0.53,1,0.1)$ & SN1 & $x=1$\\
\hline
$(0.685,1,-0.07)$ & SN2 & $x=1$\\
\hline
$(0.5,0,-0.096)$ & TR2 & $f=1/2, x=0$\\
\hline
$(0.53,0,-0.1) $& SN3 & $x=0$\\
\hline
$(0.685,0,0.07)$ & SN4 & $x=0$\\ 
\hline
\end{tabular}}
\label{tab_bifot}
\end{table}

When the human influence is strong enough, the system reaches an balanced equilibrium with 50\% of grass and 50\% of forest. On the other hand when the human influence is weak, then multiple stable steady states may occur and thus catastrophic transitions are possible between steady states with a very different percent of grass and forest as we discuss in the next section.

\subsection{Catastrophic shifts}

In an ecosystem, catastrophic shifts occur when the system passes from a stable regime to another completely different regime in an abrupt way as a consequence of environmental or external disturbances. When this happens, the recovery of the previous state is extremely difficult due (in the majority of the cases) to hysteresis effects. From a dynamical point of view, a catastrophic shift may occur in two ways: (a) as a consequence of perturbations of the vector state in a parameter region where the system has multistability, and/or (b) as a consequence of a perturbation of the value of a parameter near criticality.\\
In the first case, if the system has more than one stable regime for a specific value of the bifurcation parameter (here, the parameter $h$), then as a consequence of a perturbation in the vector state (here $x$ or $f$), the system may jump from one equilibrium to another. As each stable regime is characterized by its own basin of attraction (the set of all the initial conditions leading to the same stable regime), a perturbation to the vector state may bring it to another basin of attraction which drives the system to a different stable state.

\begin{figure}[ht]
\centerline{\includegraphics[width=5 in]{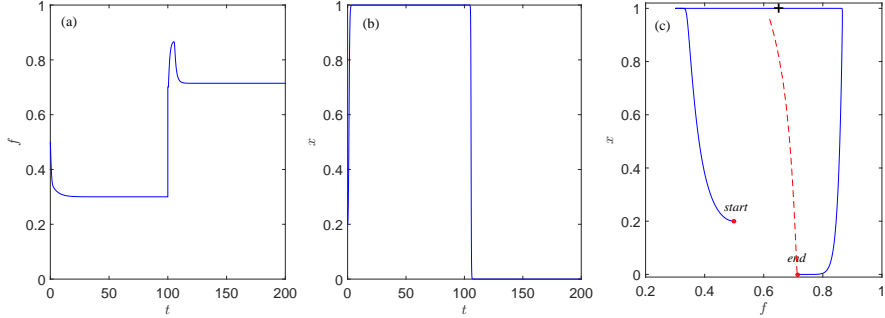}}
\vspace*{8pt}
\caption{Catastrophic shifts after perturbation of the state vector from $(f,x)=(0.33,1)\rightarrow (0.73,1)$,  at time $t=100$. Perturbation of the (a) forest state $f$ variable, (b) the grass $g$ variable. (c) The corresponding phase diagram. The dashed line is the separatrix between the attraction basins of the two stable states.}
\label{pertVec}
\end{figure}

This phenomenon is illustrated in Fig.\ref{pertVec}, where temporal simulation results are reported for an $h$ parameter value belonging to the region of multistability $[0;SN4]$. In this region, there are two stable steady states: one corresponding to $x=0$ and one to $x=1$. Starting from an initial condition within the basin of attraction of $x=1$, the system reaches first the steady state $x=1$ and then, when a perturbation is imposed (around t=100), the system finally reaches the steady state $x=0$ (see Fig.\ref{pertVec}(a),(b)). Indeed, as a consequence of the state vector perturbations, the system trajectory crosses the separatrix of the attraction basins of the two stable states (see Fig.\ref{pertVec}(c)).\\ 

\begin{figure}[ht]
\centerline{\includegraphics[width=5 in]{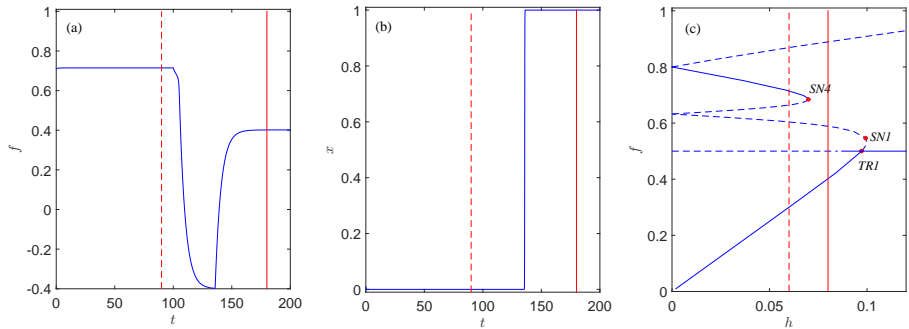}}
\vspace*{8pt}
\caption{Catastrophic shifts after perturbation of the value of $h$ (a) and (b) illustrate temporal simulations of the jump from the $x=0$ steady state to the $x=1$ steady state, respectively. (c) The bifurcation diagram of $f$ with respect to $h$. The dashed line is the value of the parameter $h$ before the perturbation and the solid line depicts the value of $h$ after the perturbation.}
\label{perth}
\end{figure}

The second way with which a catastrophic transition may occur is when one equilibrium disappears as a consequence of a perturbation of the value of the bifurcation parameter near the critical point. In this case, the system reaches a new stable equilibrium in an abrupt way. Such transitions are observed for example when a parameter crosses a critical value corresponding to a catastrophic bifurcation such as a saddle-node bifurcation point. This case is illustrated in Fig.\ref{perth}, where it is clearly shown that if the system is at the stable steady state $x=0$ just before the critical value of $SN4$, then a small perturbation of the parameter leads the system to another equilibrium $(x=1)$ with a sharp transition.

\subsection{Andronov-Hopf and Bogdanov-Takens Bifurcations} 

First, we examine the existence of periodic solutions due to the appearance of Andronov-Hopf bifurcations. Since the two branches $x=0$, $x=1$ have always real eigenvalues, we seek for Andronov-Hopf points on the third branch which is the $f=1/2$. The Jacobian matrix along the $f=1/2$ branch of solutions reads: 

\begin{equation}
J(f=1/2)=
    \quad
\begin{pmatrix} 
w'\left(1/2\right) 1/4-\nu &  \hspace{0.5cm} 2h \\
-2sx(1-x) & \hspace{0.5cm}0 
\end{pmatrix}
\label{Jac12hop}
\end{equation}
with characteristic polynomial:
\begin{eqnarray}
    P(\lambda)=|A-\lambda I| & =(w'\left(1/2\right) 1/4-\nu-\lambda)(-\lambda)+4sx(1-x) \nonumber \\ 
                             & =\lambda^2-(w'\left(1/2\right) 1/4-\nu)\lambda+4sxh(1-x)
                             \label{Charhopf}
\end{eqnarray}

If the discriminant of the characteristic polynomial is negative, i.e. $\Delta<0$ and the roots are passing the $y-$axis with non zero slope, then the assumptions for the Andronov-Hopf bifurcation are satisfied. In that case, the eigenvalues are summing up to zero or equivalently: 

\begin{equation}
    \lambda_1+\lambda_2 =Tr(J)=0 \Leftrightarrow w'\left(1/2\right) 1/4-\nu=0
    \label{casEqa}
\end{equation}
or after short calculations
\begin{equation}
    \frac{c\cdot k e^{b - k}}{\left({1 + {e^{ b- k}}}\right)^2}-\nu =0
    \label{casEqb}
\end{equation}

which is the opposite of Eq.(\ref{casNotEq}). The characteristic polynomial given by Eq.(\ref{Charhopf}) takes the form:

\begin{equation}
     P(\lambda) =\lambda^2+4sxh(1-x)
\end{equation}
with roots $\lambda_{1,2}=\pm 2 i\sqrt{sxh(1-x)}$
while $x$ must satisfy Eq.(\ref{f_12x}) and that $x h(1-x)> 0$ (taking that $s>0$). Eq.(\ref{casEqb}) is independent from $h$ and together with Eq.(\ref{f_12x}) define a system of equations which has always a solution for each $h$ with only restriction the positiveness of the product of $x h (1-x)$. Taking into account Eq.(\ref{f_12x}) we then obtain:
\begin{equation}
   xh(1-x)> 0 \Leftrightarrow 1/2 \left(1- \frac{h^*}{h} \right)\cdot h \cdot  1/2 \left(1+ \frac{h^*}{h}\right)>0 \Leftrightarrow h\left(1- \left(\frac{h^*}{h} \right)^2\right)>0 
   \label{hperiodic}
\end{equation}

In the case where $h>0$, the last inequality holds if $h>-h^{*}=h^{**}$. In order to better understand the above behaviour, we performed a two-parameter bifurcation analysis with respect to the parameters $h$, $k$.

\subsubsection{Two-parameter Bifurcation Analysis for the Andronov-Hopf Points}

In the $(h,k)$ parameter plane, the locus of the Andronov-Hopf points defined by Eq.(\ref{casEqb}) are horizontal lines since Eq.(\ref{casEqb}) is independent from $h$. With constant $\nu$, Eq.(\ref{casEqb}) is written as ($h$ satisfies Eq.(\ref{hperiodic})):

\begin{equation}
    \frac{c\cdot k e^{b - k}}{\left({1 + {e^{ b- k}}}\right)^2}-\nu =0 \Leftrightarrow \frac{z}{(1+z)^2}=\frac{\nu}{ck}
    \label{caseEq2part}
\end{equation}

where $z=e^{b-k}>0$. The function $\frac{z}{(1+z)^2}$ has two monotonic intervals with range $(0,1/4]$; then if $0<\nu/ck<1/4$ the Eq.(\ref{caseEq2part}) has two uneven roots with respect to $k$ which can be expressed in an implicit form as:

\begin{equation}
e^{b-k}=\frac{kc/\nu-2\pm \sqrt{\left(kc/\nu\right)^2-4ck/\nu}}{2}
\end{equation}
\begin{figure}[ht]
\centerline{\includegraphics[width=5 in]{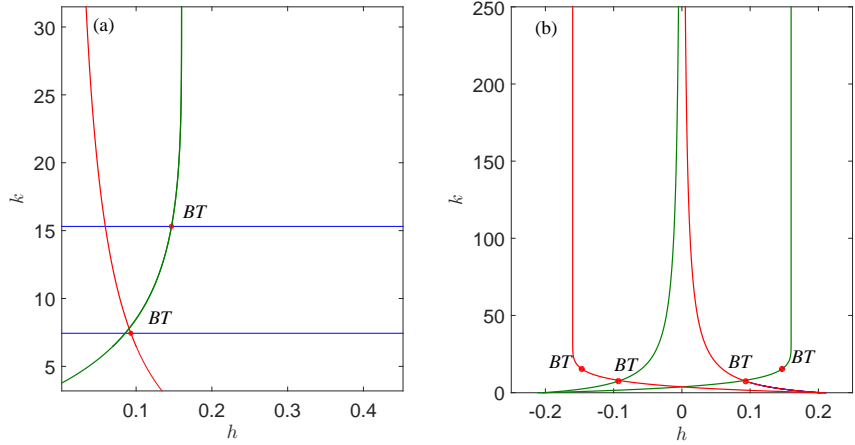}}
\vspace*{5pt}
\caption{2D bifurcation diagram with respect to $(h,k)$. \textbf{(a)} Blue horizontal lines depict the locus of Andronov-Hopf points, while the two other lines depict the locus of saddle-node bifurcations (turning points) for the $x=0$ (green) and $x=1$ (red) branch. \textbf{(b)} Continuation of limit points with respect to $h$, $k$.}
\label{2dhk_plusSaddle}
\end{figure}

For the construction of the bifurcation diagrams shown in Fig.\ref{2dhk_plusSaddle}, we used the same values for the parameters $c$, $b$, $s$, $\nu$ resulting to $k_1=7.439$ and $k_2=15.31$ (these values correspond to the horizontal lines in Fig.\ref{2dhk_plusSaddle}(a)). In the same figure, we also depict the two-parameter continuation of turning points, resulting from the system of equations Eq.(\ref{x_0}),(\ref{J110b}) and Eq.(\ref{x_1}),(\ref{J110b}). In Fig.\ref{2dhk_plusSaddle}(b), we depict the continuation of turning points for higher values of $k$ and also for $h<0$. In Fig.\ref{2dhk_plusSaddle}(b), two symmetric branches with respect to the vertical line $h=0$ exist, approaching asymptotically $k$-infinity as $h\rightarrow 0$ (see Fig.\ref{2dhk_plusSaddle}(b)).\\
An explanation of this asymptotic behavior is provided by using elements of singular two-scale analysis for the function (\ref{J110b}). The domain of forest mass variable $f$ is $f\in [0,1]$. Assuming that $k>>0$, then in the scale of order of $1/k$ i.e. $O(1/k)$, a bell shape front appears. Since $J_{11}(0)\approx -\nu$ and  $J_{11}\left( b/(b+k)\right)=\frac{c b}{4}+\frac{k-b}{2(k+b)}-\nu \rightarrow \frac{c b}{4}+\frac{1}{2}-\nu >0$ then one root of $J_{11}$ (the turning point) belongs to the interval $[0,b/(b+k)]$ i.e. $0 \leq \rho_{1}\leq b/b+k $. As $k$ increases, this implies $b/b+k \rightarrow 0 \implies \rho_{1} \rightarrow 0$. From Eq. (\ref{x_0}) or Eq.(\ref{x_1}), we obtain $h\rightarrow 0$, which is the asymptotic behavior in the neighborhood of $h=0$.\\
The second asymptotic behavior results in the slow scale ($f>1/k$). Then, the large value of $k$ vanishes the first part of $J_{11}$ concluding $J_{11}\rightarrow (1-2f)-\nu$. Setting $J_{11}=0$ in order to take the second root $\rho_{2}$, we obtain $\rho_{2}\approx \frac{1-\nu}{2}$. The corresponding value of $h$ is given from Eq.(\ref{x_0}) and Eq.(\ref{x_1}) resulting to the asymptotic values of $h$ i.e.
\begin{equation}
    h=\pm \left(c \cdot \frac{1-\nu^2}{4}- \frac{1-\nu}{2}\cdot \nu\right) 
\end{equation}
where the $\pm$ correspond to $x=0, x=1$ branches, respectively. 

Of particular interest is the case when Eq.(\ref{casEqb})  holds true along the $x=0$ or $x=1$ branches, respectively (opposite signs in (\ref{casNotEq})). Based on our analysis, this case occurs when one of the saddle-node points (on $x=1, x=0$ branches) collapse with the transcritical points at the $f=1/2$ branch (which happens when the first root of $J_{11}$,  $\rho_1=1/2$). In that case, the Jacobian matrix  has the form: 
\begin{equation}
J=
    \quad
\begin{pmatrix} 
0 &  \hspace{0.5cm} 2h \\
0 & \hspace{0.5cm} 0 
\end{pmatrix}
\label{Jac12x01all}
\end{equation}

and immediately we obtain $\lambda = 0$ as a double eigenvalue of the matrix. Thus, the critical point $(f_c, x_c, h_c, k_c)$ for the double zero eigenvalue at the $x=0$ branch can be computed from the system:
\begin{equation}
    \begin{cases}
       f=1/2\\
       x=1\\
       h=h^{*}\\
      \frac{c\cdot k e^{b - k}}{\left({1 + {e^{ b- k}}}\right)^2}-\nu =0\\
    \end{cases}       
\end{equation}

For the particular choice of the values of the other parameters, the above system of equations has the solution $(f_{c1}, x_{c1}, h_{c1}, k_{c1})=(1/2, 0, 0.146, 15.311)$.\\ Similarly for $x=1$, the system is the same, except that $x=1$ and $h=h^{**}$, with  $(f_{c2}, x_{c2}, h_{c2}, k_{c2})=(1/2, 0, 0.146, 15.311)$.\\
Setting $f \rightarrow f-1/2$, the Taylor expansion of the vector field of Eq.(\ref{system}) around the point $(1/2,0)$ reads:

\begin{eqnarray}
\left\{ \begin{array}{l}
\frac{{df}}{{dt}} = a_2(k)f^2+a_1(k)f+2hx+a_0(h,k)+h.o.t\\
\frac{{dx}}{{dt}} = -2sxf+h.o.t
\end{array} \right.
\label{systemBT}
\end{eqnarray}
where $a_2(k)=c(((16k^2e^{2(b - k)})/(e^{(b - k)} + 1)^2 + (e^{(b - k)}8(k - k^2))/4(e^{(b - k)} + 1))^2- 1/(e^{b - k)} + 1))$, $a_1(k)=((cke^{b - k})/(e^{b - k} + 1)^2 - \nu)$ and finally $a_0(k)=h - \nu/2 +\ c/4(e^{b - k} + 1)$.\\

In the 2D-parametric space $(h,k)$, the following scenarios appear:
For $h>h^{**}$ and nearby the parameter values, a horizontal line starts giving  the locus of Andronov-Hopf bifurcations i.e. Eq.(\ref{caseEq2part}). This also implies $a_1=a_0=0$, i.e. two equilibria collapse via a saddle-node bifurcation. The above cases are characteristic  of Bogdanov-Takens bifurcations (whose existence is confirmed also through the two-parameter numerical bifurcation analysis that is presented in the next section).

\section{Numerical Bifurcation Analysis}

As we have shown analytically, the simple model of forest-grassland dynamics with human interaction proposed by Innes et al.\cite{Innes_2013} exhibits a reach nonlinear behaviour with saddle-node bifurcations, transcritical bifurcations,  Andronov-Hopf bifurcations that give rise to sustained oscillations as well as codimension-two bifurcations including cusp points and Bogdanov-Takens bifurcations.

\begin{figure}[ht]
\centerline{\includegraphics[width=5 in]{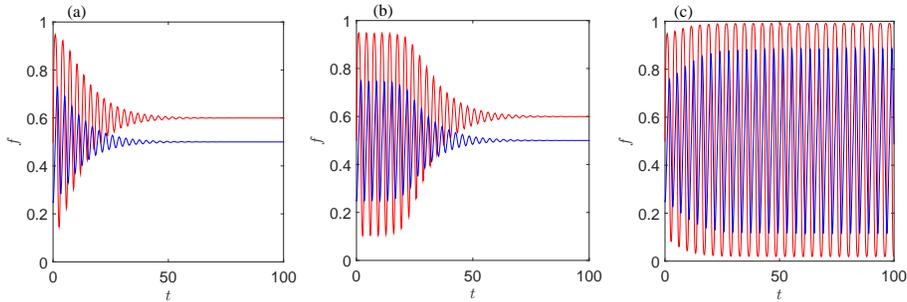}}
\vspace*{8pt}
\caption{Temporal simulations for different values of the parameter $k$, with the same initial conditions $(f,x)=(0.65,0.2)$; the other parameter values were set as $c=1, s=10, b=11$ and $h=1/2$. The blue line corresponds to $f$ while the red line to $x$  (a) $k=4$, (b) $k=4.57$ and (c) $k=4.7$.}
\label{time}
\end{figure}

With an appropriate choice of the parameter values (for example for $h=0.5$ there are transient oscillations and periodic regimes as the parameter $k$ increases. In particular, Fig.\ref{time} reports several simulations for different increasing values of $k$ (with the same initial conditions). Starting from very small values of $k$, the system quickly reaches a stable state (see Fig.\ref{time}(a)). Long transient oscillations appear as $k$ is further increased (see Fig.\ref{time}(b)). The time length of the transient oscillations also increases with the $k$ parameter and for larger values of $k$, the system manifests a permanent oscillating behavior (see Fig.\ref{time}(c)).

Although the above analytical results provide a complete picture of various types of bifurcations of fixed points, it is not possible to find other types of ``far-from-equilibrium" bifurcations such as limit (turning) points of limit cycles that can mark the onset of another type of catastrophic shift. For that purpose, one has to resort to the tools of numerical bifurcation analysis and in particular to the numerical continuation of limit cycles. This involves the solution of a boundary value problem. Let us consider for example the following system of nonlinear ordinary differential equations (that can also result from the discretization of a system of Partial Differential Equations, by e.g. finite differences or finite elements methods):

\begin{equation}\label{eqn:eq3.1}
\frac{{d\boldsymbol{x}}}{dt}\equiv \dot{\boldsymbol{x}}= \boldsymbol{F}(\boldsymbol{x}, p),\ \ \ \boldsymbol{F}: {\mathbb{R}}^n \times \mathbb{R}\rightarrow {\mathbb{R}}^n
\end{equation}

where $x \in \mathbb{R}^n$ denotes the state vector, $p\in \mathbb{R}$ is the bifurcation parameter and $F \in C^r$.

Then, one seeks solutions that satisfy the periodicity condition:

\begin{equation}\label{eqn:eq3.5}
\boldsymbol{x}(t) = \boldsymbol{x}(t + T)
\end{equation}

where $T$ denotes the (unknown) period to be computed together with $\boldsymbol{x}$. This is accomplished by solving the following rescaled with $t'=\frac{t}{T}$ boundary value problem:

\begin{equation}\label{boundaryvalue}
\begin{aligned}
\dot{\boldsymbol{x}}=T\boldsymbol{F}(\boldsymbol{x},p),  \boldsymbol{x}(0) = \boldsymbol{x}(1)\\
s.t.  \emph{h}(\boldsymbol{x}, p)=0
\end{aligned}
\end{equation}

\par
\noindent
where $\emph{h}(\boldsymbol{x}, p) = 0$ is the so-called pinning condition which factors out the infinite periodic solutions due to transnational invariance in time. 
A general pinning condition that is implemented also in MATCONT (a numerical continuation toolbox) is realized by the integral:

\begin{equation}\label{eqn:eqphasematcont}
\emph{h}(\boldsymbol{x}, p)=\int_{0}^{T}\langle \boldsymbol{x}(t), \dot{\boldsymbol{x}}_{old}(t) \rangle\,dt =0
\end{equation}

where $\dot{\boldsymbol{x}}_{old}(t)$ is the tangent vector of a previously calculated limit cycle. For the solution of the above boundary value problem, MATCONT uses orthogonal collocation, while for the continuation paste limit points, the Moore-Penrose technique is used \cite{dhooge2008new}. 

The results on the numerical bifurcation analysis of limit cycles are presented in the next section.

\subsection{Numerical-Continuation of Limit Cycles}

\begin{figure}[ht]
\centerline{\includegraphics[width=4.5 in]{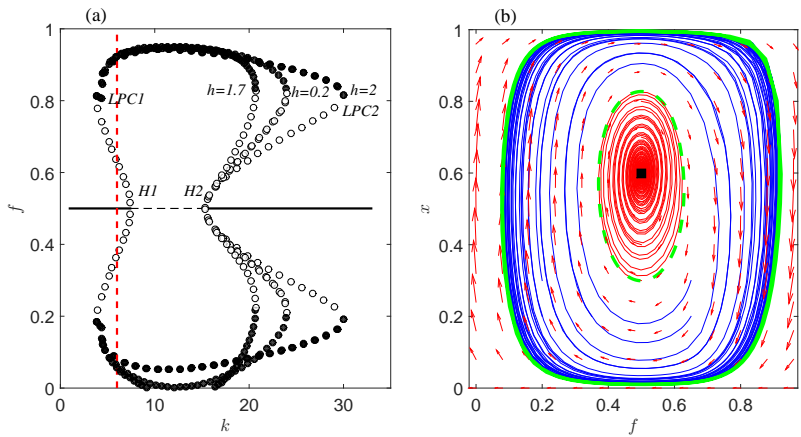}}
\vspace*{8pt}
\caption{(a) Bifurcation diagram with respect to the parameter $k$, at $\nu=0.2$, $b=11$, $c=1$ for different values of $h$. Solid lines represent stable steady states; dashed lines unstable steady states; filled circles denote stable limit cycles; empty circles denote unstable limit cycles. The points $H1$ and $H2$ represent Andronov-Hopf bifurcation points. $LPC$ represent limit point bifurcations of limit cycle. (b) Phase portrait at $k=6.5$.}
\label{bifk}
\end{figure}

We first computed numerically the one-parameter bifurcation diagram of the limit cycles by fixing different values of $h$ using $k$ as the bifurcation parameter.
We used Matcont \cite{dhooge2008new} to construct the bifurcation diagram of limit cycles with respect to the parameter $k$, for various values of $h$, $\nu=0.2$. The continuation of limit cycles with MATCONT was performed using $20$ mesh points and $4$ collocation points, with absolute and relative convergence tolerances set equal to $1E^{-04}$, thus guarantying the reliability of the numerical solution of the boundary value problem. The one-parameter bifurcation diagram is shown in Fig.\ref{bifk}. Solid lines depict stable steady solutions, while dashed lines the unstable ones. Limit cycles are indicated with circles: filled for stable and empty for unstable.

Solid lines refer to stable steady solutions, whereas dashed lines refer to unstable ones. Limit cycles are indicated with circles: filled for stable and empty for unstable.

Starting from very small values of $k$, the fixed point branch of solutions loses its stability through a subcritical Andronov-Hopf bifurcation point ($H1$) as predicted analytically. From $H1$, a branch of unstable limit cycles arises which gains stability through a turning (limit) point cycle bifurcation ($LPC1$). As noted, this point is a far-from-equilibrium bifurcation and cannot be computed without the aid of numerical bifurcation analysis. The limit cycles are stable in the range [$LPC1$; $LPC2$], whereas a branch of unstable limit cycles connects the bifurcation $LPC2$ with a second subcritical Andronov-Hopf bifurcation ($H2$). It is clear that, for this value of the parameter, there are two ranges of bistability, [$LPC1$; $H1$] and [$H2$; $LPC2$]. It is worth noting, as also shown analytically, that an increase of the value of $h$ from $0.7$ to $15$, the critical values of the $k$ parameter, where the two Andronov-Hopf bifurcations $H1$ and $H2$ occur do not change, thus leaving unchanged the range where the limit cycles are the only stable regimes. On the other hand, the critical values of $k$ where the limit point bifurcation of limit cycles $LCP2$ change, thus making the bistability windows [$H2$; $LPC2$] wider and wider, whereas the bistability range [$LPC1$; $H1$] remains the same.

\subsection{Catastrophic Shifts on Periodic Regimes}

\begin{figure}[ht]
\centerline{\includegraphics[width=5 in]{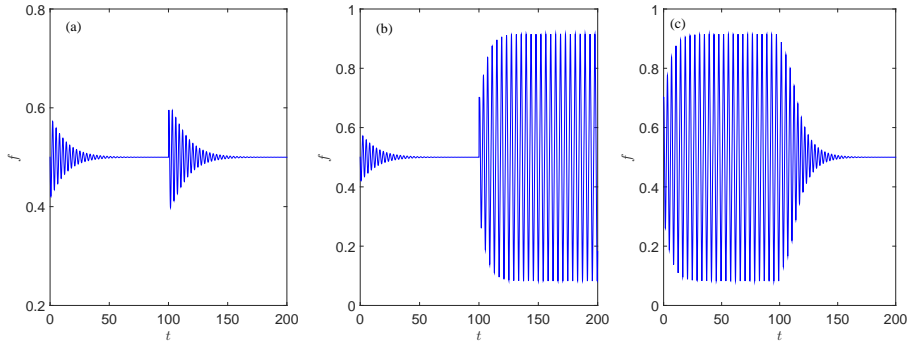}}
\vspace*{8pt}
\caption{System responses under different disturbances at $t=100$. (a) Weak perturbation of the vector space. (b) Strong perturbation of the vector space leading to an oscillation. (c) Perturbation of $k$ from $k=6.5$ to $k=3$, leading the system from an oscillatory behaviour to a steady state.}
\label{pertk}
\end{figure}

Similar analysis for the catastrophic shifts we obtain on the periodic regime. The one-parameter bifurcation diagram of limit cycles elucidates the mechanism. Fig.\ref{pertk} shows different responses of the system under different disturbances. In particular, Fig.\ref{pertk}(a), (b) show two temporal simulations after a weak perturbation (Fig.\ref{pertk}(a)) and a strong perturbation (Fig.\ref{pertk}(b)) of the state vector. In the first case, the system reaches again the final steady state, whereas with a strong perturbation, the system approaches a dynamic periodic regime. Finally, in Fig.\ref{pertk}(c), the system exhibits a catastrophic shift from the periodic regime to the steady state as consequence of a parameter perturbation.

\subsection{Two-Parameter Numerical Bifurcation Analysis}

\begin{figure}[ht]
\centerline{\includegraphics[width=4.5 in]{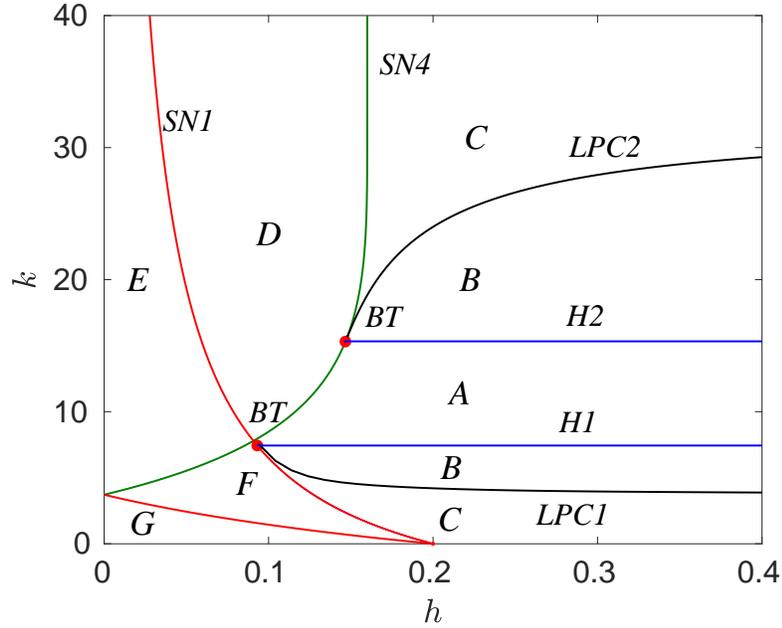}}
\vspace*{8pt}
\caption{Two-parameter bifurcation analysis with respect to the parameters $h$, $k$; the values of the other parameters are set as $b=11$, $c=1$, $s=10$, $\nu=0$. As presented in the previous sections the continuation of the turning points ($SN1$, $SN4$) and the Andronov-Hopf points was made analytically (section 2). Continuation of the limit points of the limit cycles was achieved using numerical bifurcation techniques. $BT$ denote the Bogdanov-Takens points.}
\label{2dhktot}
\end{figure}

To get a complete picture of the system's dynamics and the mechanisms that govern them, we performed a two-parameter numerical bifurcation analysis in the ($k$, $h$) parameter plane. The resulting bifurcation diagram is shown in Fig.~\ref{2dhktot}. There, we show the loci of the saddle-node bifurcations $SN1$ and $SN4$ (shown in Fig.\ref{bifx0x1}), the Andronov-Hopf bifurcations $H1$ and $H2$ and the limit point bifurcations of the limit cycles ($LPC1$) and ($LPC2$). 

It is clear that, for large values of the human impact, as reflected by the parameter $h$, the bifurcation values $LCP1$ and $LCP2$ remain constant. The loci of the various bifurcations delimit the parameter regions in different dynamics. In the case of strong human impact i.e. for higher values of $h$, the dynamics are characterized from the following scenarios (regions denoted as A, B, C in Fig.\ref{2dhktot}):
In region A we have stable limit cycles and unstable steady states. In region B we have stable and unstable limit cycles and stable  steady states. In region C we have stable steady states. We observe that as the value of $h$ increases, the corresponding critical point $LPC2$ remains almost unchanged with respect to $k,k \approx 30$, concluding that the  human interactions do not have any impact on the fire effect.
Instead, in the regime of low values of $h$, (reflecting weak human impact), the dynamics can be categorized in the following regions (see Fig.\ref{bifx0x1}). In the region D we have stable  steady states ($x=0$ branch, monostability). In the region F, we have stable  steady states ($x=1$ branch alternates with $f=1/2$ at $BT$ point of $H2$) while, the region E is characterized by multistability of steady states (see Fig.\ref{bifx0x1}). 
In the region G, we have stable  steady states ($x=1$ branch, monostability).

\section{Discussion} 
We performed an analytical and numerical one- and two-parameter bifurcation analysis of a forest grassland model which considers dynamically human interactions and their preference based on the rarity perception value. We show that despite the simplicity of the model, the nonlinear dynamic behavior is very rich including multistablity of states due to saddle-node and transcritical bifurcations, the appearance and disappearance of sustained oscillations due to Andronov-Hopf bifurcations and due to the existence of limit points of limit cycles. In this context, we derived analytically the conditions that mark the onset of such transitions and derived the loci of the bifurcations in the two-parameter space (reflecting the impact of the human influence and the effect of the fire on the rate of the reproduction of the forest). The conditions that meet codimension-two bifurcations such as Bogdanov-Takens which mark the coincide of saddle-node bifurcation (turning) points with Andronov-Hopf and saddle homoclinic bifurcations and cusp points, which mark the appearance/dissapearence of saddle-node bifurcations, are also given.  
Moreover, we performed a numerical bifurcation analysis to trace far-from-equilibrium criticalities, such as limit points of limit cycles, whose analytical computation is prohibitive. By doing so, we constructed the full two-parameter bifurcation diagram providing a full characterization of the dynamical regimes.

\bibliographystyle{ws-m3as}
\bibliography{references}

\end{document}